\documentclass[10pt]{article}
\usepackage{amsmath,amssymb,amsthm,pb-diagram,lamsarrow,pb-lams, hyperref, euscript, ulem, mathcomp}
\usepackage[matrix,arrow,curve]{xy}
\usepackage{graphicx}
\usepackage{tabularx}
\usepackage{float}
\usepackage{hyperref}

\usepackage[T1]{fontenc}

\usepackage[sc]{mathpazo}
\usepackage[usenames]{color}
\usepackage{colortbl}

\DeclareFontFamily{T1}{pzc}{}
\DeclareFontShape{T1}{pzc}{m}{it}{1.8 <-> pzcmi8t}{}
\DeclareMathAlphabet{\mathpzc}{T1}{pzc}{m}{it}

\title{Finite covering projections of noncommutative torus}

\theoremstyle{plain}
\newtheorem{prop}{Proposition}[section]

\newtheorem{lem}[prop]{Lemma}
\newtheorem{thm}[prop]{Theorem}

\theoremstyle{definition}
\newtheorem{defn}[prop]{Definition}
\newtheorem{empt}[prop]{}
\newtheorem{rem}[prop]{Remark}

\theoremstyle{remark}

\chardef\bslash=`\\ 

\newcommand{\Rb}{\mathbb{R}}

\newcommand{\eps}{\varepsilon}

\newcommand{\rar}{\rightarrow}

\newbox\ncintdbox \newbox\ncinttbox 
\setbox0=\hbox{$-$} \setbox2=\hbox{$\displaystyle\int$}
\setbox\ncintdbox=\hbox{\rlap{\hbox
    to \wd2{\hskip-.125em \box2\relax\hfil}}\box0\kern.1em}
\setbox0=\hbox{$\vcenter{\hrule width 4pt}$}
\setbox2=\hbox{$\textstyle\int$} \setbox\ncinttbox=\hbox{\rlap{\hbox
    to \wd2{\hskip-.175em \box2\relax\hfil}}\box0\kern.1em}

\begin{document}
\maketitle  \setlength{\parindent}{0pt}
\begin{center}
\author{Petr R. Ivankov*}
{\\
e-mail: * monster.ivankov@gmail.com \\
}
\end{center}

\vspace{1 in}

\begin{abstract}
\noindent

This article contains is concerned with noncommutative analogue of topological finitely listed covering projections. In my previous article \cite{ivankov:nccov_khom} I have already find a family of  covering projections of the noncommutative torus. This article describes all covering projections of the  noncommutative torus.

\end{abstract}
\tableofcontents

\section{Introduction. Preliminaries}
This article assumes elementary knowledge of following subjects
\begin{enumerate}
\item Algebraic topology  \cite{spanier:at}.
\item $C^*-$ algebras and operator theory \cite{murphy}.  

\end{enumerate}

Following notation is used.
\newline
\begin{tabular}{|c|c|}
\hline
Symbol & Meaning\\
\hline
$\mathrm{Aut}(A)$ & Group * - automorphisms of $C^*$  algebra $A$\\
$A^G$  & Algebra of $G$ invariants, i.e. $A^G = \{a\in A \ | \ ga=a, \forall g\in G\}$\\

$B(H)$ & Algebra of bounded operators on Hilbert space $H$\\
$\mathbb{C}$ (resp. $\mathbb{R}$)  & Field of complex (resp. real) numbers \\
$C(X)$ & $C^*$ - algebra of continuous complex valued \\
 & functions on topological space $X$\\
$\mathrm{Homeo}(X)$  & The space of homeomorphisms with compact-open topology \cite{spanier:at}.\\

$\mathcal{K}(H)$ or $\mathcal{K}$ & Algebra of compact operators on Hilbert space $H$\\

$\mathrm{Map}(X, Y)$  & The set of maps from $X$ to $Y$\\

  $\mathbb{T}^2$ & Commutative 2-torus  \\
$U(H) \subset \mathcal{B}(H) $ & Group of unitary operators on Hilbert space $H$\\
$U(A) \in A $ & Group of unitary operators of algebra $A$\\

$\mathbb{Z}$ & Ring of integers \\

$\mathbb{Z}_m$ & Ring of integers modulo $m$ \\

\hline
\end{tabular}
\newline
\newline
Let us recall some definitions from my previous article \cite{ivankov:nccov_khom}.

\begin{empt}{\it Galois extensions}.
Let $G$ be a finite group, $G$-Galois extensions can be regarded as particular case of Hopf-Galois extensions \cite{hajac:toknotes}, where Hopf algebra is a commutative algebra $C(G)$. Let $A$ be a $C^*$-algebra, let $G\subset \mathrm{Aut}(A)$ be a finite group of $*$- automorphisms. Let $_A\mathcal{M}^G$ be a category of $G$-equivariant modules. There is a pair of adjoint functors $(F,U)$ given by
\begin{equation}\label{f_functor}
F = A \otimes_{A^{G}} -: _{A^G}M \rightarrow _A\mathcal{M}^G;
\end{equation}
\begin{equation}\label{u_functor}
U = (-)^G: _A\mathcal{M}^G \rightarrow _{A^G}\mathcal{M}. 
\end{equation}

The unit and counit of the adjunction $(F, U)$ are given by
\begin{equation}\nonumber
\eta_{N} : N \rightarrow (A \otimes_{A^{G}} N)^{G}, \ \eta_{N}(n) = 1 \otimes  n;
\end{equation}
\begin{equation}\nonumber
\eps_{M} : A \otimes_{A^{G}} M^{G} \rightarrow M, \ \eps_{M}(a\otimes m) = am.
\end{equation}

Consider a map
\begin{equation}\label{can_def}
\mathrm{can}: A \otimes_{A^G} A \rightarrow \mathrm{Map}(G, A)
\end{equation}
given by
\begin{equation}\nonumber
a_1 \otimes a_2 \mapsto (g \mapsto a_1 (ga_2)), \ (a_1, a_2 \in A, \ g \in G).
\end{equation}
The $\mathrm{can}$ is a $_A\mathcal{M}^G$ morphism.
\end{empt}
\begin{thm}\label{hoph_galois_def_thm}\cite{morita_hopf_galois}
Let $A$ be an algebra, let $G$ be a finite group which acts on $A$, $(F, U)$ functors given by (\ref{f_functor}), (\ref{u_functor}). Consider the following
statements:
\begin{enumerate}
\item $(F,U)$ is a pair of inverse equivalences;
\item  $(F,U)$ is a pair of inverse equivalences and  $A \in _{A^{G}} \mathcal{M}$ is flat;
\item The  $\mathrm{can}$ is an isomorphism and $A \in_{A^{G}} \mathcal{M}$ is faithfully flat.

\end{enumerate}
These the three conditions are equivalent.
\end{thm}

\begin{defn}\label{hoph_galois_def}
If conditions of  theorem \ref{hoph_galois_def_thm} are hold, then $A$ is said to be a {\it left
faithfully flat $G$-Galois extension}.
\end{defn}
\begin{rem}
Theorem \ref{hoph_galois_def_thm} is an adapted to finite groups version of theorem from \cite{morita_hopf_galois}.
\end{rem}

In case of commutative $C^*$-algebras definition \ref{hoph_galois_def} supplies algebraic formulation of finitely listed covering projections of topological spaces. However I think that above definition is not quite good analogue of noncommutative covering projections. Noncommutative algebras contains inner automorphisms. Inner automorphisms are rather gauge transformations \cite{gross_gauge} than geometrical ones. So I think that inner automorphisms should be excluded. Importance of outer automorphisms was noted by  Miyashita \cite{miyashita:finite_outer_galois}. It is reasonably take to account outer automorphisms only. I have set more strong condition.   
\begin{defn}\label{gen_in_def}\cite{rieffel_finite_g}
Let  $A$ be $C^*$ - algebra. A *- automorphism $\alpha$ is said to be {\it generalized inner} if is obtained by conjugating with unitaries from multiplier algebra $M(A)$.
\end{defn}
\begin{defn}\label{part_in_def}\cite{rieffel_finite_g}
Let  $A$ be $C^*$ - algebra. A *- automorphism $\alpha$ is said to be {\it partly inner} if its restriction to some non-zero $\alpha$- invariant two-sided ideal is generalized inner. We call automorphism {\it purely outer} if it is not partly inner.
\end{defn}
Instead definitions \ref{gen_in_def}, \ref{part_in_def} following definitions are being used. 
\begin{defn}
Let $\alpha \in \mathrm{Aut}(A)$ be an automorphism. A representation $\rho : A\rightarrow B(H)$ is said to be {\it $\alpha$ - invariant} if a representation $\rho_{\alpha}$ given by
\begin{equation}
\rho_{\alpha}(a)= \rho(\alpha(a))
\end{equation}
is unitary equivalent to $\rho$.
\end{defn}
\begin{defn}
Automorphism $\alpha \in \mathrm{Aut}(A)$ is said to be {\it strictly outer} if for any $\alpha$- invariant representation $\rho: A \rightarrow B(H) $, automorphism $\rho_{\alpha}$ is not a generalized inner automorphism.
\end{defn}
\begin{defn}\label{nc_fin_cov_pr_defn}
Let $A$ be a $C^*$ - algebra and $G \subset \mathrm{Aut}(A)$ be a finite subgroup of * - automorphisms.
An injective * - homomorphism $f : A^G \rar A$ is said to be a {\it noncommutative finite
covering projection} (or {\it noncommutative $G$ - covering projection})  if $f$ satisfies following conditions:
\begin{enumerate}
\item $A$ is a finitely generated equivariant projective left and right $A^G$ Hilbert $C^*$-module.
\item If $\alpha \in G$ then $\alpha$ is strictly outer.
\item $f$ is a left faithfully flat $G$ - Galois extension.
\end{enumerate}
The $G$ is said to be {\it covering transformation group} of $f$. Denote by $G(B|A)$ covering transformation group of covering projection $A \rar B$.
\end{defn}

\section{Covering projections of noncommutative torus}

\begin{empt}\label{connected_case}{\it Noncommutative torus}. 
A noncommutative torus \cite{varilly:noncom} $A_{\theta}$ is $C^*$-norm completion of algebra generated by two unitary elements $u, v$ which satisfy following conditions
\begin{equation}\nonumber
uu^*=u^*u=vv^*=v^*v=1;\\
\end{equation}
\begin{equation}\nonumber
uv=e^{2\pi i \theta}vu
\end{equation}
where $\theta \in \Rb$.
If $\theta = 0$ then $A_{\theta}=A_0$ is commutative algebra of continuous functions on commutative torus $C(\mathbb{T}^2)$.
There is a trace $\tau_0$ on  $A_{\theta}$ such that $\tau_0 (\sum_{-\infty < i < \infty, -\infty < j<\infty}a_{ij}u^iv^j) = a_{00}$. $C^*$ - norm of $A_{\theta}$  is defined by following way $\|a\|=\sqrt{\tau_0 (a^*a)}$. 
\end{empt}
\begin{empt}\label{old_constr}
Let us recall construction from \cite{ivankov:nccov_khom}.
Let us consider * - homomorphism  $f: A_{\theta} \rightarrow A_{\theta'}$, where $A_{\theta'}$ is generated by unitary elements $u'$ and $v'$. Homomorphism $f$ is defined by following way:
\begin{equation}\nonumber
u \mapsto u'^m;
\end{equation}
\begin{equation}\nonumber
v \mapsto v'^n;
\end{equation}
It is clear  that
\begin{equation}\label{theta_k}
\theta' = \frac{\theta + k}{mn}; \ (k = 0,..., mn - 1).
\end{equation}
\end{empt}

\begin{lem}\cite{ivankov:nccov_khom}
If $\theta$ is an irrational number then above $*$-homomorphism $f: A_{\theta}\rightarrow A_{\theta'}$ is a noncommutative covering projection.
\end{lem}

\begin{empt}\label{unique_path_lifting} {\it Unique path lifting}. It is known that any topological covering projection $p: \tilde{X} \rightarrow X$ is a fibration with unique path lifting  \cite{spanier:at}, i.e. if $\omega_1, \omega_2: [0,1] \rightarrow \widetilde{X}$ are such that $\omega_1(0)=\omega_2(0)$ and $p(\omega_1(t))=p(\omega_2(t))$ ($\forall t\in [0,1]$), then $\omega_1(t)=\omega_2(t)$ ($\forall t\in [0,1]$). From unique path lifting it follows that if $\alpha: [0,1] \rightarrow \mathrm{Homeo}(X)$ is a continuous map to the space of homeomorphisms such that $\alpha(0)=\mathrm{Id}_X$ then there is the unique continuous map $\tilde{\alpha}: [0,1] \rightarrow \mathrm{Homeo}(\tilde{X})$ such that $\tilde{\alpha}(0)=\mathrm{Id}_{\tilde{X}}$, and $p(\tilde{\alpha}(t))=\alpha(t)$ ($\forall t\in[0,1]$).
\end{empt}

\begin{defn}\label{unique_path_lifting_defn}
Let $f: A^G \rightarrow A$ be a noncommutative covering projection. We say that $f$ {\it has unique lifting} if for any continuous map $\alpha: [0,1] \rightarrow \mathrm{Aut}(A^G)$ such that $\alpha(0)=\mathrm{Id}_{A^G}$ there is a map $\tilde{\alpha}: [0,1] \rightarrow \mathrm{Aut}(A)$ such that $\tilde{\alpha}|_{A^G}(t)=\alpha(t)$ ($\forall t \in [0,1]$) and $\tilde{\alpha}(0)=\mathrm{Id}_A$.
\end{defn}

\begin{empt}{\it Action of a commutative torus}.
Any point of commutative torus $\mathbb{T}^2$ can be parametrized by a pair $(z_1, z_2)\in \mathbb{C}^2$ such that $|z_1|=|z_2|=1$. Commutative torus acts on $A_{\theta}$ by following way
\begin{equation}
u \mapsto z_1u; \ v \mapsto z_2 v;  \ \forall (z_1, z_2) \in \mathbb{T}^2. 
\end{equation}
Let $G$ be a finite group and $f:A_{\theta}  \rightarrow B$ be a $G$ - covering projection, suppose that $f$ has unique lifting. Then for any $(z_1, z_2) \in \mathbb{T}^2$ there is $\alpha \in \mathrm{Aut}(B)$ such that $\alpha(a) = (z_1, z_2)a$ ($\forall a\in A_{\theta}$). Let $G'=\{\alpha\in \mathrm{Aut}(B) \ | \ \alpha|_{A_{\theta}}\in  \mathbb{T}^2 \}$.
Then there is a following exact sequence of groups
\begin{equation}\label{nc_torus_sym_grp_sequence}
\{e\}\rar G \xrightarrow{h'} G' \xrightarrow{h} \mathbb{T}^2 \rar \{e\}.
\end{equation}

Homomorphism  $h$ is a covering projection (in topological sense) because $G$ is a finite group. Covering projections of the commmutative torus are well known and exact sequence (\ref{nc_torus_sym_grp_sequence}) can be rewritten by following way
\begin{equation}\label{nc_torus_sym_grp_sequence_mod}
\{e\}\rar G_1 \times G_2 \xrightarrow{\mathrm{pr}_1 \times h'_2} G_1 \times G'_2 \xrightarrow{h} \mathbb{T}^2 \rar \{e\}
\end{equation}
where $G=G_1 \times G_2$, $G'=G_1 \times G'_2$, $G'_2$ is an abelian group which is isomorphic to $\mathbb{T}^2$, a homomorphism $G_1 \rightarrow \mathbb{T}^2 \ $ ($g \mapsto h((g, e))$) is trivial, a homomorphism $G'_2 \rightarrow \mathbb{T}^2 \ $ ($g \mapsto h((e, g ))$) is a connected covering projection of commutative torus. Sequence (\ref{nc_torus_sym_grp_sequence_mod}) can be decomposed into following sequences
\begin{equation}\label{nc_torus_sym_grp_sequence_con}
\{e\} \rar G_2 \xrightarrow{ h'_2} G'_2 \xrightarrow{h} \mathbb{T}^2 \rar \{e\};
\end{equation}
\begin{equation}\label{nc_torus_sym_grp_sequence_discon}
\{e\}\rar G_1  \xrightarrow{} G_1 \times \mathbb{T}^2 \rar \mathbb{T}^2 \rar \{e\}.
\end{equation}
Any covering of noncommutative torus can be decomposed into two covering projections which correspond to (\ref{nc_torus_sym_grp_sequence_con}) and (\ref{nc_torus_sym_grp_sequence_discon}) respectively.

Let us consider following special cases of sequence (\ref{nc_torus_sym_grp_sequence}):
\begin{enumerate}
\item (\ref{nc_torus_sym_grp_sequence_con}) $G'$ is a connected topological space $G'\rightarrow \mathbb{T}^2$ is finitely listed covering projection and $G$ is a covering transformation group.
\item (\ref{nc_torus_sym_grp_sequence_discon}) $G' = G \times  \mathbb{T}^2,$
\end{enumerate}

 \end{empt}
 
 \begin{empt}{\it 
 $G'$ is a connected topological space}.
 \newline
 In this case $G' \approx \mathbb{T}^2$. Homomorphism $h: G' \rar \mathbb{T}^2$ from (\ref{nc_torus_sym_grp_sequence}) is given by
 \begin{equation}\nonumber
 (z_1, z_2) \rar (z_1^n, z_2^m)
 \end{equation}
 where $(z_1, z_2) \in  G'\approx\mathbb{T}^2$, ($n,m \in \mathbb{N}$). Any element $(z_1, z_2) \in G \subset G'$ is given by
 \begin{equation}\nonumber
(z_1, z_2)= \left(e^{\frac{2\pi ik_1}{m}} \ e^{\frac{2\pi ik_2}{n}}\right); \ (k_1, k_2 \in \mathbb{Z}).
 \end{equation}
 Action of $G'\approx \mathbb{T}^2$ on $B$ is an unitary representation of a compact Lie group $G'\approx \mathbb{T}^2 \rightarrow \mathrm{Aut}(B)$ \cite{brocker_dieck}. Representations of $\mathbb{T}^2$ are well known, if element $b \in B$ belongs to an irreducible representation then there are $r, s \in \mathbb{Z}$ such that
  \begin{equation}\label{irr_torus}
  (z_1, z_2)b= z_1^r z_1^sb; \ ((z_1, z_2)\in G' \approx \mathbb{T}^2).
  \end{equation}
 
 An element  $a \in B$ is said to be {\it $(r, s)$ homogeneous} if it satisfies (\ref{irr_torus}).
 Let $a \in B$ be a nonzero $(r, s)$ homogeneous element, then $c=a^*a > 0$ is a $(0, 0)$ homogeneous positive element, $a$ is invariant with respect to $G$, i.e. $a \in B^G = A_{\theta}$. Any $(0, 0)$ homogeneous element in $A_{\theta}$ is a constant, i.e. $c \in \mathbb{C}$. Moreover $c \in \mathbb{R}_+$ because $c$ is positive element of $C^*$-algebra. So any $(r, s)$ homogeneous element $a$ satisfies following equation
 \begin{equation}\label{hom_repr}
 aa^* = c; \ ( \ c \in \mathbb{R}_+).
 \end{equation}
 From (\ref{hom_repr}) it follows that $\sqrt{c}a$ is an unitary, i.e any homogeneous element is $\mathbb{C}$ - proportional to an unitary element. 
 Let  $a_1, a_2 \in B$ be two unitary $(r,s)$ homogeneous elements  then $c =a_1a_2^{-1}$ is a $(0,0)$ homogeneous  element, i.e. $c\in \mathbb{C}$, or
 \begin{equation}\label{lin_dep}
 a_1 = c a_2; \ (c \in \mathbb{C});
 \end{equation}
 From (\ref{lin_dep}) that for any $(r,s)\in \mathbb{Z}^2$ a set of $(r,s)$ homogeneous elements is a one dimensional vector space over $\mathbb{C}$.
 From exactness  of $G'$ action it follows that there exist a (1, 0) homogeneous nonzero element $u'\in B$. Element $u'^m$ is $G$ invariant and $(m, 0)$ homogeneous. Element $u\in A_{\theta}$ is a  $(m, 0)$ homogeneous in $B$, so we have $u'^m=cu \ $ ($c \in \mathbb{C}$). Similarly there is an element $v'$ such that $v'^n = cv$. Monomials $u'^rv'^s$ are $(r,s)$ homogeneous elements and they are $\mathbb{C}$ - generators of $B$. From this fact it follows that any $b\in B$ can be uniquely represented by following way
 \begin{equation}
 b = \sum_{i =0; \ j=0}^{i=m-1;\ j=n-1} a_{ij}u'^iv'^j; \ (a_{ij} \in A_{\theta}).
 \end{equation}
Algebra $B$ is in fact an algebra $A_{\theta'}$ described in \ref{old_constr}.
 \end{empt}

\begin{empt} $G' = G \times  \mathbb{T}^2 $.
\newline
In this case we have following
\begin{equation}\nonumber
 G' \approx  \bigoplus_{g \in G} \mathbb{T}^2_g
\end{equation}
and a homomorphism $h: \bigoplus_{g \in G} \mathbb{T}^2_g \rightarrow \mathbb{T}^2$ is given by
\begin{equation}\nonumber
h((t_{g_1}, ..., t_{g_n})) = t_{g_1} + ... + t_{g_n}; \ (t_{g_1}, ..., t_{g_n}) \in \bigoplus_{g \in G} \mathbb{T}^2_g
\end{equation}
where additive notation of the binary $\mathbb{T}^2$ group operation is used.
The $\bigoplus_{g \in G} \mathbb{T}^2_g$ is a compact Lie group and any its representation is a direct sum of irreducible representations. Any irreducible representation $\bigoplus_{g \in G} \mathbb{T}^2_g \rar U(\mathbb{C})$ is given by
\begin{equation}\nonumber
((z_{1g_1}, z_{2g_1}), ..., (z_{1g_n}, z_{2g_n})) \mapsto z_{1g_1}^{i_{g_1}} z_{2g_1}^{j_{g_1}} \ ... \ z_{1g_n}^{i_{g_n}} z_{2g_n}^{j_{g_n}};  
\end{equation}
\begin{equation}\nonumber
(z_{1g_k}, z_{2g_k}) \in \mathbb{T}_{g_k}, \ i_{g_k}, j_{g_k} \in \mathbb{Z}.
\end{equation}
An element $a \in B$ is said to be a {\it $x = ((i_{g_1}, j_{g_1}), ..., (i_{g_n}, j_{g_n}))$ homogeneous} if it satisfies following condition
\begin{equation}\nonumber
((z_{1g_1}, z_{2g_1}), ..., (z_{1g_n}, z_{2g_n})) a = z_{1g_1}^{i_{g_1}} z_{2g_1}^{j_{g_1}} \ ... \ z_{1g_n}^{i_{g_n}} z_{2g_n}^{j_{g_n}} a.
\end{equation}
If $a'$ (resp. $a''$) is a $((i'_{g_1}, j'_{g_1}), ..., (i'_{g_n}, j'_{g_n}))$, (resp.  $((i''_{g_1}, j''_{g_1}), ..., (i''_{g_n}, j''_{g_n}))$) homogeneous element then the product $a'a''$ is a $((i'_{g_1}+i''_{g_1},j'_{g_1}+j''_{g_1}), ...,(i'_{g_n}+i''_{g_n},j'_{g_n}+j''_{g_n}))$ homogeneous element. So $B$ is a $\left(\mathbb{Z}^2\right)^G$ graded algebra. $G$ naturally acts on  $\left(\mathbb{Z}^2\right)^G$. If $x \in \left(\mathbb{Z}^2\right)^G$ and $a \in B$ is $x$ - homogeneous element then $ga$ is a $gx$ - homogeneous element. Similarly $\mathbb{T}^2$ acts on $A_{\theta}$. From this action it follows that $A_{\theta}$ is a $\mathbb{Z}^2$ graded algebra an element  $a \in A_{\theta}$ is said to be {\it $(r,s)$ homogeneous}  if
\begin{equation}\nonumber
(z_1, z_2)a=z_1^r z_2^sa.
\end{equation}
From exactness of  $\bigoplus_{g \in G} \mathbb{T}^2_g$ action it follows that there is a nonzero $((1, 0), (0, 0), ..., (0, 0))$ homogeneous element $u_{g1} \in B$. Denote by $u_g$ a homogeneous element given by
\begin{equation}\nonumber
u_g = g'u_{g_1}, \ g'g_1 = g \in G.
\end{equation}

 There is the $\mathbb{C}$ - linear map $p: B \rar A_{\theta}$ given by:
\begin{equation}\nonumber
p(a) = \sum_{g \in G}ga; \ \forall a \in B.
\end{equation}
It is clear that $p(u_{g1}) \in  A_{\theta}$ is a $(1, 0)$ homogeneous element. However any $(1, 0)$  homogeneous element is equal to $cu$ ($c\in \mathbb{C}$). If we replace $u_{g_1}$ with $c^{-1}u_{g_1}$ then $p(u_{g_1})=u$. From $p(u_{g1})p(u^*_{g1})=uu^*=1$ it follows that
\begin{equation}\label{uu_homog_product}
(u_{g_1} + ... + u_{g_n})(u^*_{g_1} + ... + u^*_{g_n}) = 1.
\end{equation}
Right part of (\ref{uu_homog_product}) is a $((0, 0), ..., (0, 0))$ homogeneous element in $B$. If $u_{g_1} u^*_{g_2} \neq 0$ then left part of (\ref{uu_homog_product}) contains a nonzero $((1, 0), (0, -1), ..., (0, 0))$ homogeneous summand but right part could not contain it, so we have $u_{g_1} u^*_{g_2} = 0$.  Similarly we can define elements $v_{g_1},..., v_{g_n}$ and  
\begin{equation}\nonumber
v_{g1} + ... + v_{g_n}= v;
\end{equation}
\begin{equation}\nonumber
(v_{g1} + ... + v_{g_n})(v^*_{g1} + ... + v^*_{g_n}) = 1.
\end{equation}
If $u_{g_1}v_{g_2} \neq 0$ then right part of
\begin{equation}\label{uv_homog_product}
uv = (u_{g1} + ... + u_{g_n})(v_{g1} + ... + v_{g_n})
\end{equation}
contains a nonzero $((1, 0), (0, 1), (0, 0), ..., (0, 0))$ homogeneous summand. However left part of (\ref{uv_homog_product}) could not contain this summand, so we have $u_{g_1}v_{g_2} = 0$.
Similarly if $g', g'' \in G$ and $g' \neq g''$ we have following:
\begin{equation}\label{zerosuuvv_product}
u_{g'}u_{g''}=u_{g'}u^*_{g''}=u^*_{g'}u_{g''}=u^*_{g'}u^*_{g''}=v_{g'}v_{g''}=v_{g'}v^*_{g''}=v^*_{g'}v_{g''}=v^*_{g'}v^*_{g''}=0;
\end{equation}
\begin{equation}\label{zerosuvvu_product}
u_{g'}v_{g''}=u_{g'}v^*_{g''}=u^*_{g'}v_{g''}=u^*_{g'}v^*_{g''}=v_{g'}u_{g''}=v_{g'}u^*_{g''}=v^*_{g'}u_{g''}=v^*_{g'}u^*_{g''}=0.
\end{equation}

From 
\begin{equation}\nonumber
u = u_{g_1} + ... +  u_{g_n}
\end{equation}
it follows that
\begin{equation}\label{w_sum}
u_{g_1}u^*_{g_1}u = w_1 + ... + w_n
\end{equation}
where $w_1$ is a $((1, 0), (0, 0), ..., (0, 0))$ homogeneous element, $w_2$ is a $((0, 0), (1, 0), ..., (0, 0))$ homogeneous element, and so on.  However from
 (\ref{zerosuuvv_product}), (\ref{zerosuvvu_product}) it follows that  $((0, 0), (1, 0), ..., (0, 0))$,...,  $((0, 0), (0, 0), ..., (1, 0))$ homogeneous summands of (\ref{w_sum}) are equal to zero, so we have
\begin{equation}\nonumber
u_{g_1}u^*_{g_1}u= w_1 = u_{g_1}
\end{equation}
or
\begin{equation}\nonumber
e_{g_1}u= u_{g_1}
\end{equation}
where $e_{g_1}=u_{g_1}u^*_{g_1}$. Similarly we can define $e_g$ for any $g \in G$ such that
\begin{equation}\label{idemp_action}
e_{g_1g_2} = g_1e_{g_2}.
\end{equation}

 From previous equations it follows that $e_g$ is an idempotent for any $g\in G$ and $B$ is a following direct sum of algebras
\begin{equation}\nonumber
B = \bigoplus_{g\in G}e_gB.
\end{equation}
A direct summand $e_gB\subset B$ is a generated by $u_g, v_g$ subalgebra. From previous   equations it follows that
 \begin{equation}\label{equ_sub}
u_gv_g=e^{2\pi i \theta}v_gu_g.
 \end{equation}
 From (\ref{equ_sub}) it follows that there is an isomorphism $A_{\theta}\rightarrow e_gB$ for any $g \in G$.
In result a noncommutative covering projection $f$ is a $*$ - homomorphism given by
 \begin{equation}\label{dir_sum}
A_{\theta} \rightarrow \bigoplus_{|G|} A_{\theta};
 \end{equation}
 \begin{equation}\nonumber
a \mapsto (a, ..., a).
 \end{equation}
From (\ref{idemp_action}) it follows that $G$ just transposes direct summands of (\ref{dir_sum}).

\end{empt}

\section{Problems}
This construction requires condition of unique path lifting \ref{unique_path_lifting_defn}. However I do not know is this condition really necessary. Analogue of infinitely listed coverings of noncommutative torus is described in \cite{ivankov:uni_nctorus}. I am engaged with the general construction of infinitely listed noncommutative covering projections.

\section{Acknowlegment}
I would like to acknowledge a "Non-commutative geometry and topology" seminar, organized by:
\begin{enumerate}
\item Prof. Alexander Mishchenko,
\item Prof. Ivan Babenko,
\item Prof. Evgenij Troitsky,
\item Prof. Vladimir Manuilov,
\item Dr. Anvar Irmatov
\end{enumerate}
for discussion of my work.   
Especially I would like acknowledge Prof. Vladimir Manuilov, because he inspired me to prove this result.

\end{document}